\documentclass[12pt]{amsart}
\usepackage{amssymb}
\usepackage{amsmath}
\usepackage{bbm}
\usepackage[dvips]{graphicx}
\usepackage{psfrag}
\usepackage{hyperref}

\numberwithin{equation}{section}

\newenvironment{proofof}[1]{\medskip\noindent
   \textbf{Proof of #1:} }{\hfill $\blacksquare$\par\medskip}

\newcommand{\be}{\begin{equation}}
\newcommand{\ee}{\end{equation}}

\newtheorem{thm}{Theorem}
\newtheorem{cor}[thm]{Corollary}
\newtheorem{lem}[thm]{Lemma}
\newtheorem{defn}{Definition}
\newtheorem{rem}{Remark}

\makeatletter
\def\Ddots{\mathinner{\mkern1mu\raise\p@
\vbox{\kern7\p@\hbox{.}}\mkern2mu
\raise4\p@\hbox{.}\mkern2mu\raise7\p@\hbox{.}\mkern1mu}}
\makeatother

\title[A bijection between permutations and a family of DPPs]
{A Natural Bijection between Permutations and a Family of 
  Descending Plane Partitions}
\author{Arvind Ayyer}
\address{Arvind Ayyer\\
Institut de Physique Th\'eorique\\
IPhT, CEA Saclay, and URA 2306, CNRS\\
91191 Gif-sur-Yvette Cedex, France}
\email{arvind.ayyer@cea.fr}
\date{\today}

\begin{document}

\begin{abstract}
  We construct a direct natural bijection between descending plane
  partitions without any special part and permutations.  The
  directness is in the sense that the bijection avoids any reference
  to nonintersecting lattice paths.  The advantage of the bijection is
  that it provides an interpretation for the seemingly long list of
  conditions needed to define descending plane partitions.
  Unfortunately, the bijection does not relate the number of parts of the
  descending plane partition with the number of inversions of the
  permutation as one might have expected from the conjecture of Mills,
  Robbins and Rumsey, although  there is a simple expression for the
  number of inversions of a permutation in terms of the corresponding
  descending plane partition.
\end{abstract}

\maketitle

\section{Introduction}
Descending plane partitions were introduced by George Andrews in order
to prove the weak Macdonald conjecture \cite{e1} and are counted by
the ASM numbers.
When they were initially introduced by Andrews, the general sense was
that these 
objects were extremely artificial and designed to specifically solve
the  conjecture.\footnote{See for example, Doron Zeilberger's
  paean to Dave Robbins \cite{z2}}

Descending plane partitions were later found to have remarkable
connections to alternating sign matrices by Mills, Robbins and Rumsey
in their proof of the Macdonald conjecture \cite{mrr1} which they
refined further in many ways in a series of conjectures in
\cite{mrr2}. Subsequently, they were  also related to 
other structures in combinatorics. Many kinds of plane partitions are
in natural bijection with classes of nonintersecting lattice paths
\cite{gv} and 
Lalonde \cite{l1} has shown, in particular, that the antiautomorphism
$\tau$ of descending plane partitions defined in \cite{mrr2} has a
natural interpretation as Gessel-Viennot path duality.  Krattenthaler
\cite{kr1} has proved a bijection between descending plane partitions
and rhombus tilings of a hexagon from which an equilateral triangle
has been removed from the center.

We will prove that the number of descending plane partitions with no
special part is the same as the number of permutations by constructing
an explicit and very natural bijection between the two objects.
Unfortunately this bijection does not naturally relate permutations
with $p$ inversions to descending plane partitions with $p$ total
parts.  We hope that a generalization of these ideas will lead to a
bijection between descending plane partitions and alternating sign
matrices.

The outline of the rest of the article is as follows. We begin with
the notations and relevant known results in Section~\ref{sec:not}. We
will need a result about descending plane partitions with one row,
which we describe in Section~\ref{sec:dpp1row} and proceed to the
proof of the bijection in Section~\ref{sec:main}. We shall give
details of the other (known) bijection through lattice paths and some
other remarks  in
Section~\ref{sec:rem}.

\section{Definitions} \label{sec:not}
We begin with a series of definitions and known results about the
objects considered here. This section is present mostly to set the
notation and experts should feel free to skip it.

\begin{defn} \label{def:dpp}
A {\em descending plane partition} (DPP) is an array $a = (a_{ij})$
of positive integers defined for $j \geq i \geq 1$ that is written in
the form
\be
\begin{array}{c c c c c c c}
a_{11} & a_{12} & \cdots & \cdots & \cdots & \cdots & a_{1,\mu_1} \\
      & a_{22} & \cdots & \cdots & \cdots & a_{2,\mu_2} & \\
& & \cdots & \cdots &\cdots & &\\
& & a_{rr} & \cdots & a_{r,\mu_r}&  &
\end{array}
\ee
where,
\begin{enumerate}
\item $\mu_1 \geq \dots \geq \mu_r$,
\item $a_{i,j} \geq a_{i,j+1}$ and $a_{i,j} > a_{i+1,j}$ whenever
  both sides are defined,
\item $a_{i,i} > \mu_i-i+1$ for $i \leq i \leq r$,
\item $a_{i,i} \leq  \mu_{i-1}-i+2$ for $1 <i \leq r$. 
\end{enumerate}
\end{defn}

The second condition in the above definition means that terms are
weakly decreasing along rows and strictly decreasing along columns.
The third condition simply means that the diagonal entry is
strictly greater than the number of entries in its row, and the fourth
condition, that it is at most
the number of entries in the row above it. Note that the last two
conditions ensure that the diagonal entries are always greater than
one.

\begin{defn}
A descending plane partition of {\em order n} is a descending plane
partition all of whose entries are less than or equal to $n$.
\end{defn}

\begin{thm} (Andrews, 1979, \cite{e1}) The number of descending plane
  partitions of order $n$, $D(n)$ is given by 
\be \label{nasm}
D(n) = \prod_{k=0}^{n-1} \frac{(3k+1)!}{(n+k)!}.
\ee
\end{thm}

We now go on to discuss refined enumeration of DPPs.
\begin{defn}
An entry $a_{i,j}$ of the  descending plane partition $a$ is
called a {\em special  part} if $a_{ij} \leq j-i$.
\end{defn}
This implies that diagonal elements can never be special parts. We
have now all definitions needed for DPPs. We go on to define ASMs and
their refinements. Another important statistic for us will be the
$r(a)$, the number of rows of the DPP $a$.

\begin{defn}
A permutation $\pi$ of the letters $\{1,\dots,n\}$ has an {\em ascent} at
  position $k$ with $1 \leq k<n$, if $\pi_k < \pi_{k+1}$.
\end{defn}

The number of permutations on $n$ letters with $k$ ascents is the
Eulerian number $E(n,k)$ \cite{gkp}, which satisfies the recurrence 
$$E(n,k) = (k+1)E(n-1,k)+(n-k)E(n-1,k-1),$$
for $n \geq 0$, $0 \leq k \leq n$ with the initial condition $E(0,k) =
\delta_{k,0}$.

\begin{defn}
The  {\em non-inversion number} $I(\pi)$ of a permutation $\pi$ on $n$ letters
is the number of pairs of elements $i,j$ such that $i<j$ and $\pi_i <
\pi_j$.
\end{defn}

The non-inversion number is the number of elementary transpositions to
convert a given permutation $\pi$ to the totally descending
permutation $n(n-1) \dots 21$, as opposed to $12 \dots n$, hence the
name. There is an obvious involution on permutations which turns
inversion numbers into non-inversion numbers.

\begin{thm} \label{thm:permdpp}
There is a natural one-to-one correspondance between descending plane
partitions of order $n$ with $k$ rows and no special part, and
permutations of size $n$ with $k$ ascents.
\end{thm}

\begin{rem}
In Theorem~\ref{thm:permdpp}, $k$ varies from 0 to $n-1$. The empty DPP,
$a = \phi$ counts as a permutation with zero rows, and vacuously, with
no special part. There is also exactly one permutation with zero
ascents, namely $\pi = n(n-1) \cdots 21$. 
\end{rem}

This immediately leads to the refined count of DPPs.

\begin{cor}
The number of descending plane partitions of order $n$ with $k$ rows and no
special parts is given by the Eulerian number $E(n,k)$.
\end{cor}

\section{Descending Plane Partitions with one row} \label{sec:dpp1row}
Before we can prove the main theorem, however, we need a simpler
result. We fix notation for future use. We denote a DPP by $a=
(a_{i,j})$ and the $i$th row of the DPP by $\alpha^{(i)}$. We will also use a
different notation for permutations suited for the interest. We will
denote a permutation with $k$ ascents by $\beta^{(1)} \cdots
\beta^{(k+1)}$, where each $\beta^{(i)}$ is decreasing. When $k=1$, we
will denote the permutation as $\beta\gamma$ to avoid unnecessary
clutter of indices.
We will also use $a_i$ and $b_i$ to denote pure numbers.

\begin{lem} \label{lem:dpp1row}
There is a natural one-to-one correspondance between descending plane
partitions of order $n$ with one row $a =(a_1,\dots,a_m)$ and
permutations of size $n$ with a single ascent $\beta \gamma$.
\end{lem}

\begin{proof}
We first associate a permutation with a single ascent to a DPP with a
single row. From the basic definitions of the DPP, we know that 
$$n\geq a_1 \geq a_2 \geq \dots \geq a_m.$$ 
Since the DPP has no
special parts, we know that $a_k \geq k$ for $1 \leq k \leq m$. But we
also know that $a_1 > m$ from the third condition in the
definition, which is stronger than the previous condition for $k=1$. 

From the DPP, we construct  
$$\gamma = (a_1,a_2-1,\dots,a_m-(m-1)).$$ 
From the weak decreasing condition above, we
clearly see that 
$$n \geq \gamma_{1}>\dots>\gamma_{m}.$$ 
From the no special part
condition, $\gamma_{k} \geq 1$ for all $k$. Therefore, $\gamma$ is a strictly
decreasing sequence of elements belonging to $[n]$. We then define
$\beta=[n] \setminus {\gamma}$ also sorted in decreasing order.
From this, we
obtain the required permutation by writing it as $\beta \gamma$. Notice that
the single ascent occurs at the junction of $\beta$ and $\gamma$
because $\gamma$ contains $m$ elements, atleast one of which is
greater than $m$, forcing at least one element not in $\gamma$ less
than $m$. 
Finally, since the
maximum value of $a_1$ is $n$, the third condition in
Definition~\ref{def:dpp} forces $m<n$ and
 $\beta$ is therefore necessarily nonempty.

The inverse procedure is quite clear. A permutation with a single
ascent can be clearly uniquely decomposed into two nonempty descending lists
$\beta$ and $\gamma$ such that the first element of $\gamma$ is
greater than the 
last element of $\beta$. 
We obtain the required DPP $a=(\gamma_1,\gamma_2+1,\dots,\gamma_m+(m-1))$. This
satisfies the 
weak decrease condition since the 
elements of $\gamma$ are strictly decreasing. Since $\gamma_k \geq 1$,
we clearly 
have $a_k \geq k$. Lastly, notice
that $\gamma_1$ has to be strictly greater than $m$, because if not, then
we are forced to have $\gamma_2=m-1,\dots,\gamma_m=1$, but this would mean that
the last element of 
$\beta$ is  $m+1$ violating the condition of a single ascent. The list
$a$ thus yields a DPP with one row and no special part.

\end{proof}

We can use this to calculate the non-inversion number for such permutations.

\begin{cor} \label{cor:inv1row}
If a permutation $\pi$ is in bijection with a descending plane
partition $a=(a_1,\dots,a_m)$ with one row and no special part, then
\be
I(\pi) = \sum_{i=1}^m a_i - m^2.
\ee
\end{cor}

\begin{proof}
We use the same notation as the proof of
Lemma~\ref{lem:dpp1row}. $I(\pi)$ is simply the total number of elementary
transpositions taken 
by the elements in $\gamma$ to return to their original
position in the completely descending permutation. We start from the
rightmost entry in $\gamma$. Clearly $\gamma_m$
will take $\gamma_m-1$ steps, $\gamma_{m-1}$ will take
$\gamma_{m-1}-2$ steps and so on. Thus 
\be
\begin{split}
I(\pi) &= (a_m-(m-1)-1)+(a_{m-1}-(m-2)-2)+ \dots + (a_1-m) \\
 &= \sum_{i=1}^m (a_i - m),
\end{split}
\ee
which gives the desired result.
\end{proof}

Notice that $\gamma$ and therefore $I(\pi)$  is determined
independently of the order of the DPP. 
For example, suppose the DPP is $a=(6,4,3)$. Then
$\gamma=(6,3,1)$. we obtain the permutation $7542631$ if $n=7$. However,
the  non-inversion number for the permutation is four, whereas $a$ has
three total parts. 
For the reverse process, consider the
permutation $25431$, which has the latter decreasing part $\gamma=5431$,
from which we get the DPP $(5,5,5,4)$.
We will need some properties of the bijection in
Lemma~\ref{lem:dpp1row} for proving Theorem~\ref{thm:permdpp}. 

\begin{lem} \label{lem:prop1row}
Using the same notation as Lemma~\ref{lem:dpp1row} and assuming $a$
has length $m$, the following hold:

\begin{enumerate}

\item $\beta_{n-m}=1$ occurs if and only if  $a_m>m$.
Assuming $1<p<n$,
$$\beta_{n-m} = p \Leftrightarrow \forall i>m-p+1, a_i=m \text{ and }
a_{m-p+1}>m.$$

\item $\beta_1 = n$ occurs if and only if $a_1<n$.
Assuming $0<p<m$,
$$\beta_{1} = n-p \Leftrightarrow \forall i \leq p,
  a_i=n \text{ and } a_{p+1}<n.$$
Lastly, $\beta_1=n-m$ if and only if $a_1 = \cdots = a_m=n$.

\end{enumerate}
\end{lem}

\begin{proof}
\begin{enumerate}
\item $\beta_{n-m}=p$ iff the letters $1,\dots,p-1$
  belong to $\gamma$, and since $\gamma$ is arranged in descending order,
  $$a_m-(m-1)=1, \cdots, a_{m-(p-2)}-(m-(p-1))=p-1,$$ 
  and furthermore $a_{m-(p-1)} -(m-p) > p$, which is precisely the
  condition stated, when $p>1$. Notice that $p$ cannot take the
  value $n$ because that would violate the single ascent condition. In
  case $p=1$, we can either have $a_m-(m-1)>1$ or $m=1$. The latter
  works because, if $m=1$, $a_1>1$ in order for the permutation to
  have a single ascent.

\item $\beta_{1}=n-p$ iff the letters $n-p+1,\dots,n$ belong to $\gamma$ and
  since $\gamma$ is arranged in descending order,
  $$a_1=n,a_2-1=n-1,\dots,a_p-(p-1)=n-(p-1),$$ 
  and the reason $n-p$ does
  not belong to $\beta$ is that either $m=p$ or $m>p$ and  $a_{p+1}-p<n-p$.
  This is again exactly the  condition stated, when $p>0$. If $p=0$,
  $n$ does not belong to $\gamma$ and thus $a_1 < n$.

\end{enumerate}
\end{proof}

\section{The Main Result} \label{sec:main}
We will construct the bijection inductively on the number of rows in
the DPP. Before that, we make some remarks on the properties of DPPs,
which follow from  Definition~\ref{def:dpp} and will need a lemma which
will be the workhorse of the proof.

\begin{rem}
\begin{enumerate}
\item Any row of a DPP is, by itself, also a valid DPP.  Moreover, a row
which is part of a DPP with no special part is also a DPP with no
special part. The latter follows from the shifted position of
successive rows.

\item Removing the last row from a DPP yields another valid
  DPP. Obviously, if the original DPP had no special part, neither
  will the new one.

\end{enumerate}
\end{rem}

\begin{lem} \label{lem:interm}
Given a set S of positive integers of cardinality $n$, there exists a
natural bijection between the DPPs, $a$, with one row and no special part
whose length $m$ satisfies $m<n$ and $a_1 \leq n$, and sequences of
all the elements of $S$ with one ascent.
\end{lem}

\begin{proof}
We define a map $\phi$ from the set $S$ to $[n]$ which takes the
smallest element to 1, the next smallest to 2 and so on until it takes
the largest element to $n$. Clearly, $\phi$ is invertible. 
Using Lemma~\ref{lem:dpp1row} therefore, we obtain a bijection between
DPPs of one row and order $n$ and no special part, and the sequence of
elements of $S$ with a single ascent. Since the DPP has order $n$, we
have $a_1 \leq n$ and therefore, the length of $a$ is strictly less
than $n$.
\end{proof}

For example, suppose $S =\{11,10,6,3,2 \}$ and $a=(4,3,3)$. The
bijection from Lemma~\ref{lem:dpp1row} yields the permutation on $n=5$
letters, $53421$, which using the map $\phi$ gives the sequence $11,6,10,3,2$.

Before we go on to the proof, we take an example to illustrate the
idea. Consider the DPP with no special part
\be
\begin{array}{c c c c c}
7  & 7  & 6 & 5 & 5 \\
   & 4  & 4 & 4 & \\
   &    & 3 & 2  & 
\end{array}
\ee
of order $n=9$, say. Then we start with the permutation
$987654321$. We will now alter it by considering the DPP row-wise. In
each row, we mark two vertical lines to separate $\beta$ and $\gamma$
using the notation in Lemma~\ref{lem:dpp1row}. The rightmost is
$\gamma$ and the one in the middle is $\beta$. The leftmost part is
completely  untouched.
\be
\begin{array}{c c c}
77655 & \to & 98|53|76421 \\
\hspace{0.15cm} 444 & \to & 9853|71|642 \\
\hspace{0.35cm} 32  & \to & 985371|4|62
\end{array}
\ee
and we end up with the permutation $985371462$, which has exactly
three ascents. In lines two and three, we have used
Lemma~\ref{lem:interm} for the rightmost part in the previous line. 

\begin{proofof}{Theorem~\ref{thm:permdpp}}
We will use induction on $k$, the number of rows of the DPP.
The case $k=1$ of the induction is precisely
Lemma~\ref{lem:dpp1row}. We now assume we have constructed, in a
one-to-one way, a permutation with $k-1$ ascents from a DPP $a$ with $k-1$
rows,
$$\alpha^{(1)},\dots,\alpha^{(k-1)}.$$
Write the permutation with $k-1$ ascents as 
$$\beta^{(1)} \cdots \beta^{(k)},$$
 where each $\beta^{(j)}$ is descending and write the $k$th row of the
DPP as $\alpha^{(k)}$.
Assume that the $k-1$th row of the DPP has length $m_{k-1}$. That is, the
terms are from $a_{k-1,k-1}$ to $a_{k-1,k+m_{k-1}-2}$. Similarly,
$\alpha^{(k)}$ has
length $m_k$, $m_k \leq m_{k-1}-1$ from Definition~\ref{def:dpp}
comprising of terms $a_{k,k}$ to $a_{k,k+m_k-1}$. 

The idea is to perform the operation on $\beta^{(k)}$ and create another
ascent within it of length $m_k$ from the right, while preserving the
ascent from $\beta_{(k-1)}$, which we describe now. Let $S$ be the set
of numbers in  $\beta^{(k)}$, which has cardinality $m_{k-1}$.
$\alpha^{(k)}$ is a DPP with one row, no special part, of length less
than $m_{k-1}$ and whose
first element $a_{k,k}$ satisfies $a_{k,k} \leq m_{k-1}$. Therefore
we are in a position to use Lemma~\ref{lem:interm} and obtain a
sequence of the elements of $S$ with a single ascent, which we call
$\gamma^{(k)}$ and $\gamma^{(k+1)}$. The length of $\gamma^{(k+1)}$ is
clearly $m_k$. We claim that this procedure is invertible and by
repeated application yields the desired permutation with $k$
ascents. What follows is a check of these claims. 

It remains to show that $\gamma^{(k)}_1$ is larger than the last entry in
$\beta^{(k-1)}$. Suppose this last entry is $p \in [n-1]$. If $p=1$, we are
done. If not, let $p \in \{2,\dots,n-1\}$. Since
the rule for creating an ascent is the same as that of creating the
first one, we can use properties of the bijection for a single
row. We will need them for the row $\alpha^{(k-1)}$.
From Lemma~\ref{lem:prop1row} (1), this implies that
$$a_{k-1,k+m_{k-1}-2}=\dots = a_{k-1,k+m_{k-1}-p}=m_{k-1}$$
 and $a_{k-1,m_{k-1}-p+1}\geq m_{k-1}+1$, and moreover that the last
$p-1$ letters of $\beta^{(k)}$ are $p-1,\dots,1$.
Thus 
$$a_{k-1,k-1} \geq \dots \geq a_{k-1,k+m_{k-1}-p-1} \geq m_{k-1}+1.$$

Notice that the first $m_{k-1}-(p-1)$ letters of $\beta^{(k)}$ are
greater than 
$p$. For it to happen that $\gamma^{(k)}_1<p$, $\gamma^{(k)}_1$ must
be one of the  
last $p-1$ letters of $\beta^{(k)}$. This implies that the action of
$\alpha^{(k)}$ 
forces all the letters larger than $\gamma^{(k)}_1$ into
$\gamma^{(k+1)}$, which can only happen if 
$a_{k,k} = \dots = a_{k,k+m_{k-1}-p}=m_{k-1}$. But this would imply
$a_{k,k+m_{k-1}-p}=a_{k-1,k+m_{k-1}-p}$, which violates condition (2)
in Definition~\ref{def:dpp}. 
Therefore the first entry of $\gamma^{(k)}$ is greater than the last entry of
$\beta^{(k-1)}$. We have thus shown that each DPP with
no special entries and with
$k$ rows gives rise to a permutation with $k$ ascents.

For the reverse process, one has to read the permutation with $k$
ascents from the right by looking at the part immediately after the
$k-1$th ascent. Using Lemma~\ref{lem:interm}, one obtains the $k$th
row of the DPP with no special parts. One then is left with a
permutation with $k-1$ ascents 
and one goes on recursively. 

Everything except the columnwise descent
is clearly ensured by this procedure. 
Essentially this occurs
because of the condition that creation of a new ascent should not kill
off an earlier ascent. 
We now describe the columnwise descent in some detail. 
We use the usual notation for the permutation with $k$ ascents,
we denote the lengths of $\beta^{(k-1)},
\beta^{(k)}$ and  $\beta^{(k+1)}$ being $m_{k-2}-m_{k-1}$,
$m_{k-1}-m_{k}$ and $m_{k}$ respectively so that the 
last three rows for the DPP, denoted $\alpha^{(k-2)},\alpha^{(k-1)}$
and $\alpha^{(k)}$ have lengths 
$m_{k-2},m_{k-1}$ and $m_{k}$ in accord with the convention used before.

We will now analyze the structure of $\alpha^{(k)}$ and
$\alpha^{(k-1)}$. In particular, we will denote the maps used in
Lemma~\ref{lem:interm} as $\phi$ and $\phi'$ respectively.

We then use the modified form of Lemma~\ref{lem:prop1row}(1) to note
that \\
$\beta^{(k-1)}_{m_{k-2}-m_{k-1}}=p$ 
implies 
$$a_{k-1,k+m_{k-1}-2} = \dots = a_{k-1,k+m_{k-1}-\phi'(p)}=m_{k-1},$$
 and 
$a_{k-1,k+m_{k-1}-\phi'(p)-1} \geq m_{k-1}+1$. This is clear because
the only change in using Lemma~\ref{lem:prop1row} directly is that
relative positions are now specified using the map $\phi'$.
Similarly, $\beta^{(k)}_{1}=r$
implies using the modified form of Lemma~\ref{lem:prop1row}(2), this
time with map $\phi$,
$$a_{k,k}= \dots = a_{k,k+m_{k-1}-\phi(r)-1} = m_{k-1},$$
 and either $m_{k-1}-m_k=r$ or
$a_{k,k+m_{k-1}-\phi(r)}<m_{k-1}-1$.  
The ascent of the permutation implies $r>p$. This in turn implies
$\phi(r) \geq \phi'(p)$ because it is possible that there are no
elements between $r$ and $p$.
A violation of the descent
condition of the DPP would entail the overlapping of the parts of the
$k-1$th and $k$th rows of $a$ which equal $m_{k-1}$. This means
 $k+m_{k-1}-\phi'(p) \leq k+m_{k-1}-\phi(r)-1$ which
implies that $\phi(r) \leq \phi'(p)-1$. But this is a contradiction. Therefore a
permutation with $k$ ascents gives rise to a DPP with $k$ rows.
\end{proofof}

We can also extend the result of Corollary~\ref{cor:inv1row} to
calculate the non-inversion number for a permutation with $k$ ascents.

\begin{cor} \label{cor:invkrow}
If a permutation $\pi$ has $k$ ascents, then the non-inversion number
 is given by the corresponding descending plane partition $a$
with $k$ rows of sizes $m_1, \dots, m_k$ as
\be
I(\pi) = \sum_{i=1}^k \sum_{j=i}^{m_i+i-1} a_{i,j} - \sum_{i=1}^k
m_i^2.
\ee
\end{cor}

\begin{proof}
Since the $i$th row of $a$ has length $m_i$, the entries are written
as $a_{i,i}, \dots, a_{i,m_i+i-1}$. 

Each successive row of the DPP is going to create more non-inversions
because one shifts successively larger numbers to the right. Moreover,
the action of each row is the same independent of the previous
rows. Therefore, one obtains the same answer for each row as in
Corollary~\ref{cor:inv1row}. Thus, the required answer is the sum for
all rows.
\end{proof}

\section{Remarks} \label{sec:rem}

We should also mention that the existence of such a bijection is part
of folklore and perhaps known to many experts, although this does not
seem to have been noted explicitly anywhere. We conjecture that
combining the bijection proposed by Gessel and Viennot \cite{gv}
between permutations and non-intersecting lattice paths with Lalonde's
bijection \cite{l1} between these paths and descending plane
partitions, one can obtain an equivalent description of the bijection
proved in this article.

\section*{Acknowledgements}
It is indeed a great pleasure to thank Robert Cori for discussions,
encouragement, and for his  hospitality at the Ecole Polytechnique. We
also thank G.X. Viennot for informative discussions and the referee
for a careful reading of the manuscript.

\bibliographystyle{alpha.bst}
\bibliography{asm}

\end{document}